\documentclass[10pt,twoside]{amsart}
\usepackage{amsmath, amsthm, amscd, amsfonts, amssymb, graphicx}
\newtheorem{thm}{Theorem}[section]
\newtheorem{cor}[thm]{Corollary}
\newtheorem{lem}[thm]{Lemma}

\newtheorem{defn}[thm]{Definition}
\newtheorem{rem}[thm]{\bf{Remark}}

\numberwithin{equation}{section}

\begin{document}
\begin{center}
{\bf{Some new Fibonacci difference spaces of non-absolute type and compact operators}}
\vspace{.5cm}
\\Anupam Das$^{1}$ and  Bipan Hazarika$^{1,\ast}$ 

\vspace{.2cm}

$^{1}$Department of Mathematics, Rajiv Gandhi University, Rono Hills, Doimukh-791 112, Arunachal Pradesh, India\\


Email: anupam.das@rgu.ac.in; bh\_rgu$@$yahoo.co.in 

\end{center}

\title{}
\author{}
\thanks{{05-04-2016\\ $^{\ast}$The corresponding author}}


\begin{abstract} The aim of the paper is to introduced the  spaces $c_{0}^{\lambda}(\hat{F})$ and $c^{\lambda}(\hat{F})$ which are the BK-spaces of non-absolute type and also derive some inclusion relations. Further, we determine the  
$\alpha-,\beta-,\gamma-$duals of those spaces and also construct their bases. We also characterize some matrix classes on the spaces $c_{0}^{\lambda}(\hat{F})$ and $c^{\lambda}(\hat{F}).$ 
Here we characterize the subclasses $\mathcal{K}(X,Y)$ of compact operators where $X$ is $c_{0}^{\lambda}(\hat{F})$ or $c^{\lambda}(\hat{F})$  and $Y$ is one of the spaces $c_{0},c, l_{\infty}, l_{1}, bv$ by applying Hausdorff measure of noncompactness.
 
\vskip 0.5cm


\vskip 0.2cm

\end{abstract}

\maketitle
\pagestyle{myheadings}
\markboth{\rightline {\scriptsize Das, Hazarika}}
         {\leftline{\scriptsize }}

\maketitle

\section{ Introduction}
\normalfont

Let $\omega$ be the space of all real-valued sequences. Any vector subspace of $\omega$ is called a $\mathit{sequence\ space}.$ By $l_{\infty},c ,c_{0},$ and $l_{p} \ (1\leq p <  \infty),$ we denote the sets of all bounded, convergent, null sequences and $p-$absolutely convergent series, respectively. Also we use the convensions that $e=(1,1,...)$ and $e^{(n)}$ is the sequence whose only non-zero term is 1 in the $nth$ place for each $n\in\mathbb{N},$ where $\mathbb{N}=\left\lbrace  0,1,2,...\right\rbrace.$

Let $X$and $Y$ be two sequence spaces and $A=(a_{nk})$ be an infinite matrix of real numbers $a_{nk},$ where $n,k\in\mathbb{N}.$ We write $A=(a_{nk})$ instead of $A=(a_{nk})_{n,k=0}^{\infty}.$ Then we say that $A$ defines a matrix mapping from $X$ into $Y$ and we denote it by writing $A:X\rightarrow Y$ if for every sequence $x=(x_{k})_{k=0}^{\infty}\in X,$ the sequence $Ax=\left\lbrace A_{n}(x) \right\rbrace_{n=0}^{\infty},$ the $A$-transform of $x,$ is in $Y,$ where
\begin{equation}
 A_{n}(x)=\sum\limits_{k=0}^{\infty} a_{nk}x_{k}\ \left( n\in \mathbb{N}\right).
 \label{1}
\end{equation}

For simplicity in notation, here and in what follows, the summation without limits runs from $0$ to $\infty.$ Also, if $x\in\omega$,then we write $x=\left(x_{k} \right)_{k=0}^{\infty}.$\par
By $(X,Y),$ we denote the class of all matrices $A$ such that $A:X\rightarrow Y.$ Thus $A\in (X,Y)$ iff the series on the right-hand side of $(1.1)$ converges for each $n\in\mathbb{N} $ and every $x\in X$ and we have $Ax\in Y$ for all  $x\in X.$ 
\par The approach constructing a new sequence space by means of matrix domain has recently employed by several authors.\par
The matrix domain $X_{A}$ of an infinite matrix $A$ in a sequence space $X$ is defined by
\begin{equation}
X_{A}=\left\lbrace x=(x_{k})\in \omega: Ax \in X\right\rbrace. 
\label{2}
\end{equation} 
Kizmaz \cite{Kizmaz} introduecd the notion of the difference operator $\Delta.$ The operator  $\Delta$ denote the matrix $\Delta=(\Delta_{nk})$ defined by
\begin{equation}
\Delta_{nk}
= \left\{
        \begin{array}{ll}
           (-1)^{n-k} , & n-1\leq k \leq n \\
           0 ,  & 0 \leq k < n-1\quad or\quad k>n.
        \end{array}
    \right.
\end{equation}
In the past, several authors studied matirx transformation on sequence spaces that are the matrix domain of the difference operator, or of the matrices of some classical methods of summability in different sequence spaces, for instance we refer to \cite{dashazarika1,Kiri,cp2,M.Mursaleen} and references therein. Hausdorff measure of noncompactness of linear operators given by infinite matrices in some special classes of sequences spaces studied by \cite{abdullah,cp4,mursallenpiloat,leen}.\\

Define the sequence $\left\lbrace f_{n}\right\rbrace_{n=0}^{\infty}$
of Fibonacci numbers given by the linear recurrence relations $ f_{0}=f_{1}=1$ and $ f_{n}=f_{n-1}+f_{n-2}, n\geq 2.$\\
Fibonacci numbers have many interesting properties and applications. For example, the ratio sequences of Fibonacci numbers converges to the golden ratio which is important in sciences and arts. Also, some basic properties of Fibonacci numbers are given as follows:
\begin{equation}
\lim_{n\rightarrow \infty}\frac{f_{n+1}}{f_{n}}= \frac{1+\sqrt{5}}{2}=\alpha\quad
(golden\ ratio),
\label{3}
\end{equation}
\begin{equation}
\sum_{k=0}^{n}f_{k}=f_{n+2}-1 \quad (n\in \mathbb{N}),
\label{4}
\end{equation}
\begin{equation}
\sum_{k}\frac{1}{f_{k}} \mbox{~converges~},
\label{5}
\end{equation}
\begin{equation}
f_{n-1}f_{n+1}-f_{n}^{2}=(-1)^{n+1} \quad (n \geq 1)  ~(Cassini ~formula)
\label{6}
\end{equation}
Substituting for $ f_{n+1} $ in Cassini's formula yields $ f_{n-1}^{2}+f_{n}f_{n-1}-f_{n}^{2}=(-1)^{n+1}.$ For details see \cite{Koshy}.

\par A sequence space $X$ is called a $FK-space$ if it is complete linear metric space with continuous coordinates $ p_{n}:X\rightarrow \mathbb{R} (n\in  \mathbb{N}),$ where $ \mathbb{R}$ denotes the real field and $ p_{n}(x)=x_{n} $  for all $ x=(x_{k}) \in X $ and every $ n \in \mathbb{N}.$ A $ BK space $ is a normed $ FK space,$ that is a $ BK-space $ is a Banach space with continuous coordinates. The sapce $ l_{p} (1\leq p < \infty)$ is a BK-sapce with\\ 
$ \parallel x \parallel_{p} =(\sum\limits_{k=0}^{\infty}\mid x_{k} \mid ^{p})^{1/p} $ \\
and $ c_{0},c $ and $ l_{\infty} $ are BK-spaces with\\ $ \parallel x \parallel_{\infty}
=\sup_{k} \mid x_{k} \mid.$\\

\par A sequence $(b_{n})$ in a normed space $X$ is called a $ Schauder \ basis $ for $X$ if every $x\in  X,$ there is a unique sequence $(\alpha_{n})$
of scalars such that $x=\sum_{n}\alpha_{n}b_{n},$ i.e., $\lim\limits_{m \rightarrow \infty}\parallel x-\sum \limits_{n=0}^{m}\alpha_{n}b_{n}\parallel =0.$\\

 The  $\alpha-,\beta-,\gamma-$duals of the sequence space $X$ are respectively defined by \\
$X^{\alpha}=\left\lbrace a=(a_{k})\in \omega : ax=(a_{k}x_{k})\in l_{1}\:\forall\: x=(x_{k})\in X  \right\rbrace. $\\
$X^{\beta}=\left\lbrace a=(a_{k})\in \omega : ax=(a_{k}x_{k})\in cs \:\forall\: x=(x_{k})\in X  \right\rbrace, $\\
and \\
$X^{\gamma}=\left\lbrace a=(a_{k})\in \omega : ax=(a_{k}x_{k})\in bs \:\forall\: x=(x_{k})\in X  \right\rbrace, $ \\
where $cs$ and $bs$ are the sequence spaces of all convergent and bounded series, respectively (See \cite{A,PK,M}).\\

\par If $X \supset \phi$ is a $BK$ space and $a \in \omega$ we write $\parallel a \parallel_{X}^{*}=\sup\left\lbrace 
\left| \sum\limits_{k=0}^{\infty} a_{k} x_{k}\right|  : \parallel x \parallel=1\right\rbrace .$\\

\par Let $X$ and $Y$ are Banach spaces.
A linear operator $L : X \rightarrow Y$ is called compact if its domain is all of $X$ and for every bounded sequence $(x_{n})_{k=0}^{\infty}$ in $X,$ the sequence $\left( L(x_{n})\right)_{n=0}^{\infty} $ has a convergent subsequence in $Y.$ We denote the class of such operators by $\mathcal{K}(X,Y).$
\par Let us recall some definitions and well-known results. 
\begin{defn}
Let $(X,d)$ be a metric space, $Q$ be a bounded subset of $X$ and $B(x,r)=\left\lbrace y\in X : d(x,y)<r\right\rbrace. $ Then the Hausdorff measure of noncompactness of $Q,$ denoted by $\chi(Q),$
is defined by
\[\chi(Q)=\inf \left\lbrace \epsilon>0: Q \subset \bigcup_{i=1}^{n} B(x_{i},r_{i}),x_{i}\in X,r_{i}< \epsilon \quad (i=1,2,...,n),n\in\mathbb{N}  \right\rbrace.\]
\end{defn}
Then the following results can be found in \cite{Bana, cp3}. 
\par If $Q,Q_{1}$ and $Q_{2}$ are bounded subsets of the metric space $(X,d),$ then we have \\
$\chi(Q)=0$ if and only if $Q$ is totally bounded set,\\
$\chi(Q)=\chi(\bar{Q}),$\\
$Q_{1}\subset Q_{2}$ implies $\chi(Q_{1})\leq \chi(Q_{2}),$\\
$\chi(Q_{1} \cup Q_{2})=$max$\left\lbrace \chi(Q_{1}),\chi(Q_{2})\right\rbrace $\\
and \\
$\chi(Q_{1} \cap Q_{2})\leq$min$\left\lbrace \chi(Q_{1}),\chi(Q_{2})\right\rbrace. $ 

\par If $Q,Q_{1}$ and $Q_{2}$ are bounded subsets of the normed space $X,$ then we have \\
 $\chi(Q_{1}+Q_{2})\leq \chi(Q_{1})+\chi(Q_{2}),$ \\
 $\chi(Q+x)=\chi(Q)$ for all $x \in X,$\\ and \\
 $\chi(\lambda Q)=\left|\lambda \right| \chi(Q)$ for all $\lambda \in \mathbb{C}.$
 \begin{defn}
 Let $X$ and $Y$ be Banach spaces and $\chi_{1}$ and $\chi_{2}$ be Hausdorff measures on $X$ and $Y.$ Then, the operator $L:X\rightarrow Y$ is called $\left( \chi_{1},\chi_{2}\right)-bounded $ if $L(Q)$ is bounded subset of $Y$ for every subset $Q$ of $X$ and there exists a positive constant $K$ such that $\chi_{2}(L(Q))\leq K \chi_{1}(Q)$ for every bounded subset $Q$ of $X.$ If an operator $L$ is $(\chi_{1},\chi_{2} )-$bounded then the number 
 $\parallel L\parallel_{(\chi_{1}, \chi_{2})}
 =\inf \left\lbrace K>0 : \chi_{2}(L(Q))\leq K \chi_{1}(Q) \mbox{~for~all~bounded~} Q \subset X\right\rbrace $ is called $(\chi_{1}, \chi_{2})-$ measure of noncompactness of $L.$ In particular, if $\chi_{1}=\chi_{2}=\chi,$ then we write 
 $\parallel L\parallel_{(\chi, \chi)}=\parallel L \parallel_{\chi}.$
 \end{defn}
 
  \par The idea of compact operators between Banach spaces is closely related to the Hausdorff measure of noncompactness, and it can be given as follows:
 
 \par Let $X$ and $Y$ be Banach spaces and $L \in B(X,Y).$ Then the Hausdorff measure of noncompactness of $L,$ is denoted by $\parallel L\parallel_{\chi},$ can be given by
 \begin{equation}
 \parallel L\parallel_{\chi}=\chi(L(S_{X}))
 \end{equation}
 where $S_{X}=\left\lbrace x\in X : \parallel x\parallel=1\right\rbrace $
 and we have  $L$ is compact if and only if
\begin{equation}
 \parallel L\parallel_{\chi}=0
 \end{equation}
 We also have\\
 $\parallel L \parallel = \sup_{x \in S_{X}}\parallel L(x)\parallel_{Y}.$

\section{The sequence spaces $c_{0}^{\lambda}(\hat{F})$ and $c^{\lambda}(\hat{F})$ of non-absolute type }

In this section, we introduce the spaces $c_{0}^{\lambda}(\hat{F})$ and $c^{\lambda}(\hat{F})$   and show that these spaces are the BK-spaces of non-absolute type which are linearly isomorphic to the spaces  $c_{0}$ and $c,$ respectively.
\par
We shall assume throughout this paper that $\lambda =\left( \lambda_{k} \right) _{k=0}^{\infty} $ is strictly increasing sequence of positive reals tending to $\infty,$ that is 
$0< \lambda_{0} <\lambda_{1}<$ ... and $\lambda_{k} \rightarrow \infty$ as $k \rightarrow \infty .$
\par
Recently, the sequence spaces $c_{0}^{\lambda}$ and $c^{\lambda}$ of non absolute type have been introduced by Mursaleen and Noman (see \cite{19}) as follows:
\[c^{\lambda}_{0}=\left\lbrace x=(x_{k})\in w : \lim\limits_{n}\frac{1}{\lambda_{n}}\sum\limits_{k=0}^{n}(\lambda_{k}-\lambda_{k-1})x_{k} =0\right\rbrace\]
and
\[c^{\lambda}=\left\lbrace x=(x_{k})\in w : \lim\limits_{n}\frac{1}{\lambda_{n}}\sum\limits_{k=0}^{n}(\lambda_{k}-\lambda_{k-1})x_{k}\mbox{~exists~}\right\rbrace.\]

Also, it has been shown that the inclusions $c_{0} \subset c_{0}^{\lambda},c \subset c^{\lambda} $ and $c_{0}^{\lambda} \subset c^{\lambda}$ hold.\\

Let $f_{n}$ be the $nth$ Fibonacci number for every $n \in \mathbb{N}$. The infinite matrix  $\hat{F}=\left( f_{nk} \right) $
was introduced by Kara \cite{Kara} is defined as follows.
\begin{equation}
f_{nk}
= \left\{
        \begin{array}{ll}
           -\frac{f_{n+1}}{f_{n}} , 
           & k=n-1 \\
            \frac{f_{n}}{f_{n+1}} ,
           & k=n \\
           0,
           & 0\leq k < n-1\quad or\quad k>n
        \end{array}
    \right.
\end{equation}
where $n,k \in \mathbb{N}.$ 

\par Define the sequence $y=(y_{n}),$ which will be frequently used, by the $\hat{F}$-transform of a sequence $x=(x_{n}),$ i.e., $y_{n}=\hat{F}_{n}(x),$ where
\begin{equation}
y_{n}
= \left\{
        \begin{array}{ll}
           \frac{f_{0}}{f_{1}}x_{0}=x_{0} , 
           &  n=0 \\
            \dfrac{_{f_{n}}}{f_{n+1}}x_{n}-\dfrac{_{f_{n+1}}}{f_{n}}x_{n-1} ,
           &  n \geq 1
        \end{array}
    \right.
  \label{11}  
\end{equation}
where $ n \in \mathbb{N}.$
\par We employ a technique of obtaining a new sequence space by means of matrix domain. We thus introduce sequence spaces $c_{0}^{\lambda}(\hat{F})$ and $c^{\lambda}(\hat{F})$ are defined as follows.
\[c_{0}^{\lambda}(\hat{F})=\left\lbrace x=(x_{k})\in \omega : \lim\limits_{n} \frac{1}{\lambda_{n}}\sum\limits_{k=0}^{n} \left( \lambda_{k}-\lambda_{k-1} \right) \left( \dfrac{{f_{k}}}{f_{k+1}}x_{k}-\dfrac{{f_{k+1}}}{f_{k}}x_{k-1}    \right)=0  \right\rbrace \]
and
\[c^{\lambda}(\hat{F})=\left\lbrace x=(x_{k})\in \omega : \lim\limits_{n} \frac{1}{\lambda_{n}}\sum\limits_{k=0}^{n} \left( \lambda_{k}-\lambda_{k-1} \right) \left( \dfrac{{f_{k}}}{f_{k+1}}x_{k}-\dfrac{{f_{k+1}}}{f_{k}}x_{k-1}    \right) \mbox{~exists~} \right\rbrace.\]

We shall use the convention that any term with negative subscript is equal to zero, e.g. $\lambda_{-1}=0$ and $x_{-1}=0.$ 
\par We can redefine the spaces $c_{0}^{\lambda}(\hat{F})$ and 
$c^{\lambda}(\hat{F})$ by
\begin{equation}\label{c}
c_{0}^{\lambda}(\hat{F})=\left( c_{0}^{\lambda}\right)_{\hat{F}} 
\mbox{~and~} 
c^{\lambda}(\hat{F})=\left( c^{\lambda}\right)_{\hat{F}}.
\end{equation}

 It is immediate by (\ref{c}) that the sets $c_{0}^{\lambda}(\hat{F})$ and $c^{\lambda}(\hat{F})$ are linear spaces with coordinatewise addition and scalar multiplication.  On the other hand, we define the matrix $\bar{F}=\left( \bar{f}_{nk} \right) $ for all $n,k \in \mathbb{N}$ by
\begin{equation}
\bar{f}_{nk}
= \left\{
      \begin{array}{ll}
                \frac{1}{\lambda_{n}}\left[ \left( \lambda_{k}-\lambda_{k-1}\right)  \frac{f_{k}}{f_{k+1}}-\left( \lambda_{k+1}-\lambda_{k}\right)  \frac{f_{k+2}}{f_{k+1}}
                \right]   
                 &   ,k < n \\
                \frac{1}{\lambda_{n}} \left( \lambda_{n}-\lambda_{n-1}\right)  \frac{f_{n}}{f_{n+1}}
                 &  ,k=n \\
                 0 & ,k>n.
              \end{array}
\right.
\end{equation}

Then, it can be easily seen that 
\begin{equation}\label{c0}
\bar{F}_{n}(x)
= \frac{1}{\lambda_{n}}\sum\limits_{k=0}^{n} \left( \lambda_{k}-\lambda_{k-1} \right) \left( \dfrac{{f_{k}}}{f_{k+1}}x_{k}-\dfrac{{f_{k+1}}}{f_{k}}x_{k-1}    \right)
\end{equation}
holds for all $n \in \mathbb{N}$ and $x=(x_{k})\in w,$ which leads us to the fact that 
\begin{equation}\label{c1}
c_{0}^{\lambda}(\hat{F})=\left( c_{0}\right)_{\bar{F}} 
\mbox{~and~}
c^{\lambda}(\hat{F})=\left( c\right)_{\bar{F}}.
\end{equation}
\par Moreover, it is obvious that $\bar{F}$ is a triangle. Thus it has a unique inverse 
 $\bar{F}^{-1}=\left( \bar{f}_{nk}^{-1} \right) $ for all $n,k \in \mathbb{N}$ given by

\begin{equation}
\bar{f}_{nk}^{-1}
= \left\{
      \begin{array}{ll}
                \lambda_{k}f_{n+1}^{2}\left[ \frac{1}{ \lambda_{k}-\lambda_{k-1}}. \frac{1}{f_{k}f_{k+1}}-\frac{1}{ \lambda_{k+1}-\lambda_{k}}. \frac{1}{f_{k+1}f_{k+2}}
                \right]   
                 &   ,0 \leq k < n \\
               \lambda_{n}f_{n+1}^{2}.\frac{1}{ \lambda_{n}-\lambda_{n-1}}. \frac{1}{f_{n}f_{n+1}}
                 &  ,k=n \\
                 0 & ,k>n.
              \end{array}
\right.
\end{equation}

\par Further, for any sequence $x=(x_{k})$ we define the sequence $y=(y_{k})$ such that $y=\bar{F}(x)$ and we observe that
\begin{equation} \label{c2}
y_{k}=\bar{F}_{k}(x)=\sum\limits_{j=0}^{k-1}
\frac{1}{\lambda_{k}} \left[ \left( \lambda_{j}-\lambda_{j-1} \right) \dfrac{{f_{j}}}{f_{j+1}}- \left( \lambda_{j+1}-\lambda_{j} \right) \dfrac{{f_{j+2}}}{f_{j+1}} \right]x_{j} 
+\frac{1}{\lambda_{k}} \left( \lambda_{k}-\lambda_{k-1} \right) \dfrac{{f_{k}}}{f_{k+1}}x_{k},
\end{equation}
where $k \in \mathbb{N}.$
\par Now, we may begin with the following theorem which is essential in the text.
\begin{thm}
The sequence spaces $c^{\lambda}_{0}(\hat{F})$ and $c^{\lambda}(\hat{F})$ are BK-spaces with norm 
\[\parallel x \parallel_{c^{\lambda}_{0}(\hat{F})}= \parallel x \parallel _{c^{\lambda}(\hat{F})}=
\parallel \bar{F}(x) \parallel_{l_{\infty}}
=\sup \limits_{n} \left| \bar{F}_{n}(x) \right|.\]
\end{thm}
\begin{proof}
Since (\ref{c1}) holds and $c_{0}$ and $c$ are BK-spaces with respect to their natural norm and the matrix $\bar{F}$ is a triangle. Theorem 4.3.12 of Wilansky \cite{cp6} gives the fact that $c^{\lambda}_{0}(\hat{F})$ and $c^{\lambda}(\hat{F})$ are BK-spaces with given norms.
\end{proof}

\begin{rem}
One can easily check that the absolute property is not satisfied by  $c^{\lambda}_{0}(\hat{F})$ and $c^{\lambda} (\hat{F}),$ that is $\parallel x \parallel _{c^{\lambda}_{0}(\hat{F})} \neq \parallel \mid x \mid  \parallel _{c^{\lambda}_{0}(\hat{F})}$
  and
$\parallel x \parallel _{c^{\lambda}(\hat{F})} \neq \parallel \mid x \mid  \parallel _{c^{\lambda}(\hat{F})}.$ This shows that 
 $c^{\lambda}_{0}(\hat{F})$ and $c^{\lambda}(\hat{F})$  are sequence spaces of non-absolute type, where $\mid x \mid = (\mid x_{k} \mid).$
\end{rem}

\begin{thm}
The sequence spaces $c^{\lambda}_{0}(\hat{F})$ and $c^{\lambda}(\hat{F})$ of non-absolute type are linearly isomorphic to the spaces $c_{0}$ and $c,$ respectively, that is 
$c^{\lambda}_{0}(\hat{F}) \cong c_{0}$ and $c^{\lambda}(\hat{F}) \cong c.$
\end{thm}
\begin{proof}
To prove this, we should show the existence of a linear bijection between the spaces  $c^{\lambda}_{0}(\hat{F})$ and $c_{0}.$ Consider the transformation $T$ defined, with the notation of (\ref{c2}), from $c^{\lambda}_{0}(\hat{F})$ to $c_{0}$ by 
$Tx =y=\bar{F}(x) \in c_{0}$ for every $x \in c^{\lambda}_{0}(\hat{F}).$ Also, the linearity of $T$ is clear. Further, it is trivial that $x=0$ whenever $Tx=0$ hence $T$ is injective.
\par Further, let $y \in (y_{k})\in c_{0}$ and we define the sequence $x=(x_{k})$ by
\begin{equation}\label{c3}
x_{k}=\sum\limits_{j=0}^{k}\sum\limits_{i=j-1}^{j} (-1)^{j-i} 
\frac{\lambda_{i}y_{i}}{\lambda_{j}-\lambda_{j-1}}.\frac{f_{k+1}^{2}}{f_{j}f_{j+1}}
\end{equation}
for $k=0,1,2,... $ and so on
and $\bar{F}_{n}(x)=y_{n}.$ This shows that $\bar{F}(x)=y$ and since $y \in c_{0},$ we obtain  $\bar{F}(x) \in c_{0}.$ Thus, we deduce that  $x \in c^{\lambda}_{0}(\hat{F})$ and $Tx=y.$ Hence $T$ is surjective.
\par Moreover, for every $x \in c^{\lambda}_{0}(\hat{F})$ we have
\[\parallel Tx\parallel_{c_{0}}=\parallel Tx\parallel_{l_{\infty}}
=\parallel y\parallel_{l_{\infty}}=\parallel \bar{F}(x)\parallel_{l_{\infty}}=\parallel x\parallel_{c^{\lambda}_{0}{(\hat{F})}}\]
which means that $T$ is norm preserving. Consequently, $T$ is a linear bijection which shows that $ c^{\lambda}_{0}(\hat{F})$ and $c_{0}$ are linearly isomorphic.
\par Similarly, we can show that 
$c^{\lambda}(\hat{F}) \cong c $  and this concludes the proof.
\end{proof}

\begin{thm}
The space $l_{\infty}$ does not include the spaces $c^{\lambda}_{0}(\hat{F})$ and $c^{\lambda}(\hat{F}).$ 
\end{thm}
\begin{proof}
We have, from  equation (\ref{c0}) that \\
$\bar{F}_{n}(x)
= \frac{1}{\lambda_{n}}\sum\limits_{k=0}^{n} \left( \lambda_{k}-\lambda_{k-1} \right) \left(\dfrac{{f_{k}}}{f_{k+1}}x_{k}-\dfrac{{f_{k+1}}} {f_{k}}x_{k-1}    \right) =\frac{1}{\lambda_{n}}\sum \limits_{k=0}^{n} \left( \lambda_{k}-\lambda_{k-1} \right)\hat{F}_{k}(x).$
\par 
Let us define a sequence $x=(x_{k})=(f_{k+1}^{2}).$ Since 
$f_{k+1}^{2} \rightarrow \infty$ as $k \rightarrow \infty$ and 
$\hat{F}(x)=e^{(0)}=(1,0,0,...),$ therefore
 $\bar{F}_{n}(x)=\frac{\lambda _{0}}{\lambda_{n}} \rightarrow 0$ as
 $n \rightarrow \infty.$
 \par 
 Hence  we can conclude that $x \in c^{\lambda}_{0}(\hat{F}) $ but not in $l_{\infty}.$ Similarly, we can show that $x \in c^{\lambda}(\hat{F}) $ but not in $l_{\infty}.$
\end{proof}

\begin{thm}
The inclusion $c^{\lambda}_{0}(\hat{F}) \subset c^{\lambda}(\hat{F})$ strictly holds.
\end{thm}
\begin{proof}
It is clear that  $c^{\lambda}_{0}(\hat{F}) \subseteq c^{\lambda}(\hat{F}).$ Consider the sequence $x=(x_{n})$ defined by
\begin{equation}
x_{n}
= \left\{
        \begin{array}{ll}
           1  
           & , n=0 \\
            f_{n+1}^{2}\left( \sum\limits_{j=1}^{n} \frac{1}{f_{j}f_{j+1}}+1\right) 
           &  ,n \geq 1
        \end{array}
    \right.
  \label{c4}  
\end{equation}
Then, we have $\bar{F}(x)=e$ and hence $\bar{F}(x) \in c \smallsetminus c_{0}$ where $e=(1,1,1,...).$ Thus the sequence $x$ in $c^{\lambda}(\hat{F})$ but not in $c^{\lambda}_{0}(\hat{F}).$ Hence the inclusion $c^{\lambda}_{0}(\hat{F}) \subset c^{\lambda}(\hat{F})$ strict.
\end{proof}
\begin{thm}
The inclusion $c \subset c^{\lambda}(\hat{F})$ and $c_{0} \subset c^{\lambda}_{0}(\hat{F})$ strictly hold.
\end{thm}
\begin{proof}
Let $x \in \left( x_{k}\right) \in c.$ We have $c \subset c^{\lambda}$ and $c \subset c_{0}^{\lambda}$ if and only if 
$\liminf\limits_{n \rightarrow \infty}\frac{\lambda_{n+1}}{\lambda_{n}}=1$ (see \cite{19}) and $\hat{F}(x)\in c $ as $\lim\limits_{n\rightarrow \infty}\hat{F}_{n}(x)$ exists, therefore,  $\hat{F}(x)\in c^{\lambda}. $ This shows that $x\in c^{\lambda}(\hat{F}).$ Consequently,  $c \subseteq c^{\lambda}(\hat{F}).$
\par Let $x=\left( x_{k} \right)=(f_{n+k}^{2})\notin c $ but we have $\hat{F}(x)=\left( 1,0,0,0,... \right) $ so $y=(y_{k}),$ where
$y_{k}=\bar{F}_{k}(x)=\frac{\lambda_{0}}{\lambda_{k}}.$ Thus $y \in c,$ hence $x\in  c^{\lambda}(\hat{F})$ but not in $c.$ We conclude that  $c \subset c^{\lambda}(\hat{F})$ hold strictly.
\par Similarly, we can show that 
$c_{0} \subset c_{0}^{\lambda}(\hat{F})$ hold strictly. 
\end{proof}

\par Now, because of the transformation $T$ defined from $c^{\lambda}_{0}(\hat{F})$ to $c_{0},$ is an isomorphism, the inverse image of the basis $\left\lbrace e^{(k)} \right\rbrace_{k=0}^{\infty} $ of the space  $c_{0}$ is the basis for the new space $c^{\lambda}_{0}(\hat{F}).$ Therefore, we have the following result.

\begin{thm}
Define the sequence $b^{(k)}=\left\lbrace b_{n}^{(k)} \right\rbrace_{n=0}^{\infty}$ for every fixed $k=0,1,2,...$ by
\begin{equation}
b_{n}^{(k)}
= \left\{
        \begin{array}{ll}
           0 
           &  ,n < k \\
            \frac{\lambda_{k}}{\lambda_{k}-\lambda_{k-1}}. \frac{f_{n+1}^{2}}{f_{k}f_{k+1}} 
           &   ,n=k \\
           \frac{\lambda_{k}}{\lambda_{k}-\lambda_{k-1}}. \frac{f_{n+1}^{2}}{f_{k}f_{k+1}}-\frac{\lambda_{k}}{\lambda_{k+1}-\lambda_{k}}. \frac{f_{n+1}^{2}}{f_{k+1}f_{k+2}}
           &  ,n>k
        \end{array}
    \right.
  \label{c5}  
\end{equation}
 where $ n=0,1,2,... $. Then the sequence $ (b^{(k)})_{k=0}^{\infty}$ is a basis for the space $ c^{\lambda}_{0}(\hat{F}),$ and every $  x \in  c^{\lambda}_{0}(\hat{F})$ has a unique representation of the form 
\begin{equation}
 x= \sum\limits_{k} \alpha_{k}b^{(k)}.
 \label{c6}
 \end{equation}
 where $\alpha_{k}=\bar{F}_{k}(x)$ for all $k=0,1,2,...$.
 \label{t1}
\end{thm}

\begin{proof}
It is clear that the inclusion $\left\lbrace b^{(k)}\right\rbrace \subset c_{0}^{\lambda}(\hat{F}) $ holds, since 
\begin{equation}
\bar{F}\left( b^{(k)} \right) = e^{(k)} \in c_{0}, \: k \in \mathbb{N}.
\label{c7}
\end{equation}
Let $x \in c_{0}^{\lambda}(\hat{F})$ be given. For every non-negative integer $ m,$ we put
$x^{(m)}= \sum\limits_{k=0}^{m} \alpha_{k}b^{(k)}.$ Then we obtain by (\ref{c7}) that
$\bar{F}\left( x^{(m)}\right)=\sum \limits_{k=0}^{m} \alpha_{k}\bar{F}\left( b^{(k)}\right) =\sum \limits_{k=0}^{m} \bar{F}_{k}(x)e^{(k)} $ \\
and hence \\
\begin{equation}
\bar{F}_{n}(x-x^{(m)})
= \left\{
        \begin{array}{ll}
           0  
           & ,0 \leq n \leq m \\
            \bar{F}_{n}(x)  ,
           &  n > m
        \end{array}
    \right. 
    \label{c8}
\end{equation}
where $n,m \in \mathbb{N}.$

\par Now, given $\epsilon > 0,$ there exists a non-negative integer $m_{0}$ such that
 $\left|\bar{F}_{m}(x) \right| < \epsilon/2$ 
for all $m \geq m_{0}.$
\par 
Therefore, we have for every $ m \geq m_{0} $ that
\[ \parallel x-x^{(m)}\parallel_{c_{0}^{\lambda}(\hat{F})} =  \sup_{n>m} \left| \bar{F}_{n}(x)\right| \leq 
\sup_{n>m_{0}} \left| \bar{F}_{n}(x)\right| \leq \frac{\epsilon}{2}
 < \epsilon.\]
which shows that
 $ \lim \limits_{m \rightarrow \infty} \parallel x-x^{(m)} \parallel_{c_{0}^{\lambda}(\hat{F})}=0$ 
  and hence $x$ is represented as in (\ref{c6}).
 \par 
 If possible let there exists another representation
 
 \begin{equation}
  x= \sum\limits_{k} \beta_{k}b^{(k)}.
  \label{c9}
  \end{equation}
  \par Since $T \equiv \bar{F}$ is a linear transformation from $c_{0}^{\lambda}(\hat{F})$ to $c_{0}$ and is continuous, therefore we have
  $\bar{F}_{n}(x)=\sum\limits_{k}\beta_{k}\bar{F}_{n}\left( b^{(k)} \right)=\beta_{n}, $ for all $n=0,1,2,...$. Thus we have
  $\alpha_{k}=\beta_{k}$ for all $k=0,1,2,....$ Hence the representation (\ref{c6}) is unique.
\end{proof}
\begin{thm}
The sequence $\left\lbrace b,b^{(0)},b^{(1)},... \right\rbrace $ is a basis for the space $c^{\lambda}(\hat{F})$ and every $x \in c^{\lambda}(\hat{F})$ has unique representation of the form,
\begin{equation}
x=lb+\sum\limits_{k}\left( \alpha_{k}-l\right)b^{(k)} 
\label{c10}
\end{equation}
where $\alpha_{k}=\bar{F}_{k}(x)$ for all $k=0,1,2,...$, the sequence $b=(b_{n})$ is defined by 
\begin{equation}
b_{n}
= \left\{
        \begin{array}{ll}
           1  
           & , n=0 \\
            f_{n+1}^{2}\left( \sum\limits_{j=1}^{n} \frac{1}{f_{j}f_{j+1}}+1\right) 
           &  ,n \geq 1,
        \end{array}
    \right.
  \label{c11}  
\end{equation}
the sequence 
$b^{(k)}=\left\lbrace b_{n}^{(k)}\right\rbrace_{n=0}^{\infty} $ is defined by (\ref{c5})
for every fixed $k=0,1,2,...$ and 
\begin{equation}
l=\lim\limits_{k}\bar{F}_{k}(x).
\label{c12}
\end{equation}
\end{thm}
\begin{proof}
Since $\left\lbrace b^{(k)}\right\rbrace \subset c_{0}^{\lambda}(\hat{F}) $ and $\bar{F}(b)=e \in c,$ the inclusion $\left\lbrace b, b^{(k)}\right\rbrace \subset  c^{\lambda}(\hat{F}) $ trivially holds. Further, let $x \in c^{\lambda}(\hat{F}).$ Then there exists a unique $l$ satisfying
(\ref{c12}). Thus  we have $y \in c_{0}^{\lambda}(\hat{F}),$ where $y=x-lb.$ Therefore, by Theorem \ref{t1} we have that the representation 
$y=\sum\limits_{k}\beta_{k}b^{(k)}$ is unique, where
$\beta_{k}=\bar{F}_{k}(x-lb)=\alpha_{k}-l$ for all $k.$ Hence the representation (\ref{c10}) is unique.
\end{proof}
\begin{cor}
The difference spaces $c_{0}^{\lambda}(\hat{F})$ and $c^{\lambda}(\hat{F})$ are separable.
\end{cor}
\section{ The $\alpha$-,$\beta$- and $\gamma$-duals of the spaces $c_{0}^{\lambda}(\hat{F})$ and $c^{\lambda}(\hat{F})$}
In this section, we determine the $\alpha-$,$\beta-$ and $\gamma-$duals of the sequence space $c_{0}^{\lambda}(\hat{F})$ and $c^{\lambda}(\hat{F})$ of non-absolute type.   \\

We shall assume throughout our discussion that  the sequences $x=(x_{k})$ and $y=(y_{k})$ are connected by the relation (\ref{c2}). Now we may begin with quoting the following lemmas (see \cite{Michael}) which are needed to prove next theorems.
\begin{lem}
$A \in \left(c_{0}:l_{1} \right)=\left(c:l_{1} \right)$ if and only if 
$$\sup_{K \in \mathcal{F}} \sum \limits_{n}\left| \sum\limits_{k \in 
K}a_{nk}\right|< \infty. $$
\label{l1}
\end{lem}

\begin{lem}
$A \in \left(c_{0}:c \right)$ if and only if
\begin{equation}
\lim \limits_{n}a_{nk} \mbox{~exists~for~each~} k \in \mathbb{N},
\label{c13}
\end{equation}
\begin{equation}
\sup_{n}\sum\limits_{k} \left| a_{nk}\right|<\infty 
\label{c14}
\end{equation}
\label{l2}
\end{lem}
\begin{lem}
$A \in \left( c:c \right)$ if and only if (\ref{c13}) and (\ref{c14}) hold, and 
\begin{equation}
\lim \limits_{n}\sum \limits_{k}a_{nk} \mbox{~exists~}.
\label{c15}
\end{equation}
\label{l3}
\end{lem}
\begin{lem}
$A \in \left( c_{0}: l_{\infty}\right)=\left(c: l_{\infty} \right)$ if and only if (\ref{c14}) holds.
\label{l4}
\end{lem}
Now, we prove the following results.
\begin{thm}
The $\alpha$-dual of the sequence space $c_{0}^{\lambda}(\hat{F})$ and $c^{\lambda}(\hat{F})$ is the set \\
$$ b_{1}=\left\lbrace a=(a_{k}) \in \omega :  \sup_{K \in \mathcal{F}} \sum \limits_{n}\left| \sum\limits_{k \in K} b_{nk}\right| < \infty
 \right\rbrace,$$  
where the matrix $B=(b_{nk})$  is defined via the sequence $a=(a_{n})$ by 
 $$ b_{nk}= \left\{
         \begin{array}{ll}
                   \left( \frac{\lambda_{k}}{\lambda_{k}-\lambda_{k-1}}. \frac{f_{n+1}^{2}}{f_{k}f_{k+1}}-\frac{\lambda_{k}}{\lambda_{k+1}-\lambda_{k}}. \frac{f_{n+1}^{2}}{f_{k+1}f_{k+2}}\right) a_{n} 
                   &  ,k < n \\
                    \frac{\lambda_{k}}{\lambda_{k}-\lambda_{k-1}}. \frac{f_{n+1}^{2}}{f_{k}f_{k+1}}.a_{n} 
                   &   ,k=n \\
                   0                    &  ,k>n
        \end{array}
  \right. $$ for all $n,k \in \mathbb{N}$
  and $a=(a_{n}) \in \omega.$
\end{thm}
\begin{proof}
Let $a=(a_{n}) \in \omega.$ Then by (\ref{c2}) and (\ref{c3}) we immediately derive that
\begin{equation}\label{c16}
a_{n}x_{n}=\sum\limits_{k=0}^{n}\sum\limits_{j=k-1}^{k} (-1)^{k-j} 
\frac{\lambda_{j}y_{j}}{\lambda_{k}-\lambda_{k-1}}.\frac{f_{n+1}^{2}}{f_{k}f_{k+1}}a_{n}=B_{n}(y),
\end{equation}
where $n=0,1,2,...$. Thus we observed that by (\ref{c16}) that $ax=\left( a_{n}x_{n}\right) \in l_{1}$ when $x=(x_{k})\in c_{0}^{\lambda}(\hat{F}) $ or $c^{\lambda}(\hat{F})$ if and only if $By \in l_{1}$ when $y=(y_{k})\in c_{0}$ or $c$ i.e. $a=(a_{n})$ is in the $\alpha-$ dual of the spaces $c_{0}^{\lambda}(\hat{F})$ or $c^{\lambda}(\hat{F})$ if and only if $B \in \left( c_{0} : l_{1}\right)=\left( c : l_{1}\right).$ We, obtain by Lemma \ref{l1} that 
$a \in \left\lbrace c_{0}^{\lambda}(\hat{F})\right\rbrace^{\alpha} =
\left\lbrace c^{\lambda}(\hat{F})\right\rbrace^{\alpha} $ iff 
$$ \sup_{K \in \mathcal{F}} \sum\limits_{n}\left| \sum\limits_{k \in 
K}b_{nk}\right|< \infty $$ which gives
$\left\lbrace c_{0}^{\lambda}(\hat{F})\right\rbrace^{\alpha} =
\left\lbrace c^{\lambda}(\hat{F})\right\rbrace^{\alpha}=b_{1}.$
\end{proof}

\begin{thm}\label{T1}
Define the sets $b_{2},b_{3},b_{4}$ and $b_{5}$ by \\
$ b_{2}=\left\lbrace a=(a_{k}) \in \omega :  \sum\limits_{j=k}^{\infty}  a_{j} f_{j+1}^{2} \mbox{~exists~for~each~} k\in \mathbb{N}   \right\rbrace,$ \\
$ b_{3}=\left\lbrace a=(a_{k}) \in \omega : \sup_{n} \sum\limits_{k=0}^{n-1}\left| \bar{a}_{k}(n)\right|< \infty \right\rbrace $,
\\ 
$ b_{4}=\left\lbrace a=(a_{k}) \in \omega : \sup_{n}\left| \frac{\lambda_{n}}{\lambda_{n}-\lambda_{n-1}}.\frac{f^{2}_{n+1}}{f_{n}f_{n+1}}.{a}_{n}\right|< \infty \right\rbrace $\\
and\\
$ b_{5}=\left\lbrace a=(a_{k}) \in \omega : a_{0}+ \sum\limits_{k=1}^{\infty}\left\lbrace f_{k+1}^{2}\left( \sum\limits_{j=1}^{k}\frac{1}{f_{j}f_{j+1}}+1\right) \mbox{~converges~} \right\rbrace 
 \right\rbrace ;$ \\
 where\\
 $\bar{a}_{k}(n)
 =\lambda_{k}\left[ \frac{a_{k}}{\lambda_{k}-\lambda_{k-1}}. \frac{f_{k+1}^{2}}{f_{k}f_{k+1}}+
 \left( \frac{1}{\lambda_{k}-\lambda_{k-1}}. \frac{1}{f_{k}f_{k+1}}-
 \frac{1}{\lambda_{k+1}-\lambda_{k}}. \frac{1}{f_{k+1}f_{k+2}}
 \right)\sum\limits_{j=k+1}^{n}f_{j+1}^{2}a_{j}
 \right] $ , $k<n.$\\
 Then $ \left\lbrace  c_{0}^{\lambda}(\hat{F})\right\rbrace ^{\beta}=b_{2}\cap b_{3}\cap b_{4} $ and 
 $ \left\lbrace  c^{\lambda}(\hat{F})\right\rbrace ^{\beta}=b_{3}\cap b_{4}\cap b_{5}. $
\end{thm}     
\begin{proof}
Let $a=(a_{k})\in \omega$ and consider the equality, \\
$
\sum\limits_{k=0}^{n} a_{k}x_{k}\\
=
\sum\limits_{k=0}^{n} \left\lbrace \sum\limits_{j=0}^{k}
\left[ \sum\limits_{i=j-1}^{j} (-1)^{j-i} 
\frac{\lambda_{i}y_{i}}{\lambda_{j}-\lambda_{j-1}}.\frac{f_{k+1}^{2}}{f_{j}f_{j+1}}\right] \right\rbrace a_{k}\\
=\sum\limits_{k=0}^{n-1} \lambda_{k}\left[ \frac{a_{k}}{\lambda_{k}-\lambda_{k-1}}.\frac{f_{k+1}^{2}}{f_{k}f_{k+1}} + \left( \frac{1}{\lambda_{k}-\lambda_{k-1}}.\frac{1}{f_{k}f_{k+1}}-
\frac{1}{\lambda_{k+1}-\lambda_{k}}.\frac{1}{f_{k+1}f_{k+2}}\right)\sum\limits_{j=k+1}^{n}f_{j+1}^{2}a_{j} \right] y_{k} +
\frac{\lambda_{n}y_{n}}{\lambda_{n}-\lambda_{n-1}}.\frac{f_{n+1}^{2}}{f_{n}f_{n+1}}a_{n}\\
=\sum\limits_{k=0}^{n-1}\bar{a}_{k}(n)y_{k}+
\frac{\lambda_{n}y_{n}}{\lambda_{n}-\lambda_{n-1}}.\frac{f_{n+1}^{2}}{f_{n}f_{n+1}}a_{n}\\
 =T_{n}(y); 
 $ where $n=0,1,2,3,...$ and $T=\left( t_{nk}\right) $ is defined by \\
  $$ t_{nk}
 = \left\{
          \begin{array}{ll}
                    \bar{a}_{k}(n) 
                    &  ,k < n \\
                     \frac{\lambda_{n}}{\lambda_{n}-\lambda_{n-1}}. \frac{f_{n+1}^{2}}{f_{n}f_{n+1}}.a_{n} 
                    &   ,k=n \\
                    0
                    &  ,k>n
         \end{array}
   \right. $$ for all $n,k \in \mathbb{N}.$\\
   
\par Then we have $ax=(a_{k}x_{k})\in cs$ whenever $x=(x_{k}) \in c_{0}^{\lambda}(\hat{F})$ iff $Ty \in c$ whenever $y=(y_{k}) \in c_{0}.$ Therefore $a=(a_{k}) \in \left\lbrace  c_{0}^{\lambda}(\hat{F})\right\rbrace ^{\beta} $ iff $T \in \left( c_{0} : c\right). $ Therefore by using Lemma \ref{l2} we derive that 
\begin{equation}\label{c17}
\sum\limits_{j=k}^{\infty} a_{j}f_{j+1}^{2}  \mbox{~exists~each~} k=0,1,2,...
\end{equation}
\begin{equation}\label{c18}
\sup_{n}\sum\limits_{k=0}^{n-1}\left| \bar{a}_{k}(n)\right| < \infty
\end{equation}
and
\begin{equation}\label{c19}
\sup_{n} \left| \frac{\lambda_{n}}{\lambda_{n}-\lambda_{n-1}}.\frac{f_{n+1}^{2}}{f_{n}f_{{n+1}}}.a_{n}\right| < \infty.
\end{equation}   
Hence we conclude that $\left\lbrace c_{0}^{\lambda}\left(\hat{F} \right) \right\rbrace^{\beta}
=b_{2}\cap b_{3}\cap b_{4}. $  \\
 
Similarly, from Lemma \ref{l3} we have $a=(a_{k}) \in \left\lbrace c^{\lambda}\left(\hat{F} \right) \right\rbrace^{\beta} $ if and only if $T \in (c : c).$ Therefore, we derive from (\ref{c13}),(\ref{c14}) that (\ref{c17}),(\ref{c18}) and (\ref{c19}) hold.
\par Further, it can easily be seen that the equality
\begin{equation}\label{c20}
a_{0}+ \sum\limits_{k=1}^{n}\left\lbrace f_{k+1}^{2}\left( \sum\limits_{j=1}^{k} \frac{1}{f_{j}f_{j+1}}+1\right)a_{k} \right\rbrace
= \sum\limits_{k=0}^{n-1}\bar{a}_{k}(n)+
\frac{\lambda_{n}}{\lambda_{n}-\lambda_{n-1}}.\frac{f_{n+1}^{2}}{f_{n}f_{n+1}}a_{n}
=\sum\limits_{k} t_{nk}
\end{equation}
where $n=0,1,2,...$.
Consequently, we obtain from Lemma \ref{l3} that
\begin{equation}\label{c21}
a_{0}+ \sum\limits_{k=1}^{n}\left\lbrace f_{k+1}^{2}\left( \sum\limits_{j=1}^{k} \frac{1}{f_{j}f_{j+1}}+1\right)a_{k} \right\rbrace \mbox{~converges~}.
\end{equation}
Thus $\sum\limits_{j=k}^{\infty}a_{j}f_{j+1}^{2}$ exists is a weaker condition than $a_{0}+ \sum\limits_{k=1}^{n}\left\lbrace f_{k+1}^{2}\left( \sum\limits_{j=1}^{k} \frac{1}{f_{j}f_{j+1}}+1\right)a_{k} \right\rbrace$
converges. Hence we conclude that 
$\left\lbrace  c^{\lambda}(\hat{F})\right\rbrace ^{\beta}=b_{3}\cap b_{4}\cap b_{5}.$
\end{proof}

\begin{thm}
The $\gamma-$dual of the spaces $ c_{0}^{\lambda}(\hat{F})$ and
$ c^{\lambda}(\hat{F})$ is the set
$b_{3}\cap b_{4}.$
\end{thm}
\begin{proof}
This theorem can be proved similarly as the proof of the Theorem \ref{T1} by using Lemma \ref{l4}.
\end{proof}
\section{ Some matrix mappings related to the sequence spaces $ c_{0}^{\lambda}(\hat{F})$ and
$ c^{\lambda}(\hat{F})$ }

In this section, we characterize the classes 
$ \left( c^{\lambda}(\hat{F}): l_{p} \right),$ 
$ \left( c_{0}^{\lambda}(\hat{F}): l_{p} \right),$
$ \left( c^{\lambda}(\hat{F}): c \right),$
$ \left( c_{0}^{\lambda}(\hat{F}): c \right),$
$ \left( c^{\lambda}(\hat{F}): c_{0} \right)$ and
$ \left( c_{0}^{\lambda}(\hat{F}): c_{0}
 \right),$
 where $1 \leq p \leq \infty .$ 
 \par We assume that the sequences $x=(x_{k})$ and $y=(y_{k})$ are connected by $y=\bar{F}(x).$ Also for an infinite matrix $A=(a_{nk}),$ we shall write that \\
 $\bar{a}_{nk}(m)
  =\lambda_{k}\left[ \frac{a_{nk}}{\lambda_{k}-\lambda_{k-1}}. \frac{f_{k+1}^{2}}{f_{k}f_{k+1}}+
  \left( \frac{1}{\lambda_{k}-\lambda_{k-1}}. \frac{1}{f_{k}f_{k+1}}-
  \frac{1}{\lambda_{k+1}-\lambda_{k}}. \frac{1}{f_{k+1}f_{k+2}}
  \right)\sum\limits_{j=k+1}^{m}f_{j+1}^{2}a_{nj}
  \right],$ \\ where  $k<m$ and\\
  $\bar{a}_{nk}
   =\lambda_{k}\left[ \frac{a_{nk}}{\lambda_{k}-\lambda_{k-1}}. \frac{f_{k+1}^{2}}{f_{k}f_{k+1}}+
   \left( \frac{1}{\lambda_{k}-\lambda_{k-1}}. \frac{1}{f_{k}f_{k+1}}-
   \frac{1}{\lambda_{k+1}-\lambda_{k}}. \frac{1}{f_{k+1}f_{k+2}}
   \right)\sum\limits_{j=k+1}^{\infty}f_{j+1}^{2}a_{nj}
   \right]$\\
 for all $n,k,m \in \mathbb{N}$ provided the convergence of the series.
 \par The following lemmas will be needed in our discussion.
  
\begin{lem} \cite{cp6} \label{l5} 
The matrix mapping between the BK-spaces are continuous.
\end{lem}
\begin{lem} \cite{Michael} \label{l6}
$A\in \left( c : l_{p}\right)$ if and only if 
$\sup_{F \in \mathcal{F}} \sum\limits_{n}\left| \sum\limits_{k \in F} a_{nk}\right|^{p} < \infty$ where $1 \leq p < \infty.$
\end{lem} 
\begin{lem} \cite{Michael}\label{l7} 
$A \in (c : c_{0})$ if and only if 
\begin{equation}
\sup_{n}\sum\limits_{k}\left| a_{nk} \right| < \infty, 
\label{c22}
\end{equation}
\begin{equation}
\lim\limits_{n}a_{nk}=0 \: for \:all \: k \in \mathbb{N},
\label{c23}
\end{equation}
\begin{equation}
\lim\limits_{n}\sum\limits_{k}a_{nk}=0. 
\label{c24}
\end{equation}
\end{lem}
\begin{lem}  \cite{Michael}\label{l8} 
$A \in \left( c_{0}, c_{0}\right) $ if and only if (\ref{c22}) and
(\ref{c23}) hold.
\end{lem}
Now, we prove the following results.
\begin{thm}\label{T2}
(i) Let $1 \leq p < \infty$. Then $A \in \left( c^{\lambda}\left( \hat{F}\right) : l_{p} \right) $ if and only if
\begin{equation}
\sup_{F \in \mathcal{F}}\sum\limits_{n}\left|\sum\limits_{k \in F}\bar{a}_{nk} \right|^{p} < \infty,
\label{c25}
\end{equation} 
\begin{equation}
\sup_{m}\sum\limits_{k=0}^{m-1}\left| \bar{a}_{nk}(m)\right| < \infty ; (n \in \mathbb{N}),
\label{c26}
\end{equation}
\begin{equation}\label{c27}
a_{n0}+ \sum\limits_{k=1}^{\infty}\left\lbrace f_{k+1}^{2}\left( \sum\limits_{j=1}^{k} \frac{1}{f_{j}f_{j+1}}+1\right)a_{nk} \right\rbrace \mbox{~converges~}, n \in \mathbb{N},
\end{equation}
\begin{equation}\label{c28}
\lim\limits_{k \rightarrow \infty} \frac{f_{k+1}^{2}}{f_{k}f_{k+1}}.\frac{\lambda_{k}}{\lambda_{k}-\lambda_{k-1}}a_{nk}=a_{n}; (n \in \mathbb{N}),
\end{equation}
\begin{equation}\label{c29}
(a_{n})\in l_{p}.
\end{equation}
(ii) $A \in \left( c^{\lambda}\left( \hat{F}\right): l_{\infty}  \right) $ if and only if (\ref{c27}) and (\ref{c28}) hold, and
\begin{equation}\label{c30}
\sup_{n}\sum\limits_{k} \left| \bar{a}_{nk}\right| < \infty,
\end{equation}
\begin{equation}\label{c31}
\left( a_{n}\right)\in l_{\infty}.
\end{equation}
\end{thm}
\begin{proof}
Suppose that the conditions from (\ref{c25}) to (\ref{c29}) hold and let $x=(x_{k})\in c^{\lambda}\left(\hat{F}\right).$ Then we have $\left\lbrace a_{nk}\right\rbrace _{k \in \mathbb{N}}\in \left\lbrace c^{\lambda}\left(\hat{F}\right)\right\rbrace ^{\beta} $ for all $n \in \mathbb{N}$ and this implies that $Ax$ exists. Also, it is clear that the associated sequence $y=(y_{k})$ such that $T(y)=x$ in the space $c$ and $y_{k}\rightarrow l$ as
$k \rightarrow \infty$ for some suitable $l.$ Combining the Lemma \ref{l6} with (\ref{c25}), we get 
$\bar{A}=(\bar{a}_{nk})\in \left( c: l_{p}\right),$ $1 \leq p < \infty.$
\par Let us now consider the following equality derived by using the relation $y=\bar{F}(x)$ and the $mth$ partial sum of the series $\sum\limits_{k}a_{nk}x_{k}$  then we have
\begin{equation}\label{c32}
\sum\limits_{k=0}^{m}a_{nk}x_{k}
=\sum\limits_{k=0}^{m-1}\bar{a}_{nk}(m)y_{k}+\frac{f_{m+1}^{2}}{f_{m}f_{m+1}}.\frac{\lambda_{m}}{\lambda_{m}-\lambda_{m-1}}.a_{nm}y_{m}
\end{equation}
where $m,n \in \mathbb{N}.$
\par Since $y \in c$ and $\bar{A}\in \left( c: l_{p}\right);$ $
\bar{A}y$ exists and so the series $\sum\limits_{k} \bar{a}_{nk}y_{k}$ converges for all $n \in \mathbb{N}.$ From (\ref{c27}) we get $\sum\limits_{j=k}^{\infty}a_{nj}$ converges for all $n,k \in \mathbb{N}$ and hence $\bar{a}_{nk}(m) \rightarrow \bar{a}_{nk}$ as $m \rightarrow \infty.$ Therefore, as $m \rightarrow \infty$ we get from (\ref{c32}) and (\ref{c28}) that
\begin{equation}\label{c33}
\sum\limits_{k}a_{nk}x_{k}
=\sum\limits_{k}\bar{a}_{nk}y_{k}+la_{n} \mbox{~for~all~} n \in \mathbb{N}
\end{equation}
 and which can be written as
\begin{equation}\label{c34}
A_{n}(x)=\bar{A}_{n}(y)+la_{n} \mbox{~for~all~} n \in \mathbb{N}.
\end{equation}
 Therefore we have
\begin{equation}\label{c35}
\parallel A(x)\parallel_{l_{p}} \leq 
\parallel \bar{A}(y)\parallel_{l_{p}}+ \mid l \mid \parallel(a_{n})\parallel_{l_{p}} < \infty
\end{equation}
which gives $Ax \in l_{p}$ and hence $A \in \left( c^{\lambda}\left( \hat{F}\right):l_{p} \right),$ where $1 \leq p < \infty.$
\par Conversely suppose that $A \in \left( c^{\lambda}\left( \hat{F}\right):l_{p} \right),$ where $1 \leq p < \infty.$ Then 
$\left\lbrace a_{nk}\right\rbrace \in \left\lbrace c^{\lambda}\left( \hat{F}\right) \right\rbrace^{\beta} $ for all $n \in \mathbb{N}$ and this with Theorem \ref{T1} implies the necessity conditions (\ref{c26}) and (\ref{c27}).
\par Since $c^{\lambda}\left( \hat{F}\right)$ and $l_{p}$ are $BK-$spaces, we have by Lemma \ref{l5} that there is a constant $M > 0$ such that 
\begin{equation}\label{c36}
\parallel Ax\parallel_{l_{p}} \leq M \parallel x\parallel_{c^{\lambda}\left( \hat{F}\right)}
\end{equation}
for all $x \in c^{\lambda}\left( \hat{F}\right).$ Now let $F \in \mathcal{F}.$ Then the sequence $z=\sum\limits_{k\in F}b^{(k)}$ is in $c^{\lambda}\left( \hat{F}\right), $ where the sequence 
$b^{(k)}=\left\lbrace b^{(k)}_{n}\right\rbrace _{n \in \mathbb{N}} $ for every fixed $k \in \mathbb{N}.$ We have by (\ref{c7}) that\\
$\parallel z \parallel_{c^{\lambda}\left( \hat{F}\right)}
=\parallel \bar{F}(z) \parallel_{l_{\infty}}
=\parallel \sum\limits_{k \in F}\bar{F}(b^{(k)}) \parallel_{l_{\infty}}
=\parallel \sum\limits_{k \in F}e^{(k)} \parallel_{l_{\infty}}=1.$
\\ Again for every $n \in \mathbb{N},$ we have 
$A_{n}(z)=\sum\limits_{k \in F}A_{n}\left( b^{(k)}\right)
=\sum\limits_{k \in F}\sum\limits_{j}a_{nj}b_{j}^{(k)}
=\sum\limits_{k \in F} \bar{a}_{nk}.$
Since the inequality (\ref{c36}) is satisfied for the sequence $z\in c^{\lambda}\left( \hat{F}\right),$ we have for any $F \in \mathcal{F}$ that
$\left( \sum\limits_{n} \left| \sum\limits_{k \in F} \bar{a}_{nk}\right|^{p} \right)^{1/p} \leq M$ which implies  $\sup_{F\in \mathcal{F}}\sum\limits_{n}\left|
\sum\limits_{k\in F} \bar{a}_{nk} \right|^{p}< \infty.$ Thus it follows by Lemma \ref{l6} that $\bar{A}=(\bar{a}_{nk})\in (c : l_{p}).$
\par Now let $y=(y_{k})\in c\smallsetminus c_{0}$ and consider the sequence $x=(x_{k})$ defined by (\ref{c3}) for every $k\in \mathbb{N}.$
Then $x \in c^{\lambda}\left( \hat{F}\right) $ such that $y=\bar{F}(x).$ Therefore $Ax$ and $\bar{A}y$ exist. This gives the series $\sum\limits_{k}a_{nk}x_{k}$ and $\sum\limits_{k}\bar{a}_{nk}y_{k}$ 
converges for all $n \in \mathbb{N}.$
\par We also have
\[\lim\limits_{m}\sum \limits_{k=0}^{m-1} \bar{a}_{nk}(m)y_{k}=\sum\limits_{k} \bar{a}_{nk}y_{k}; \: (n \in \mathbb{N}).\]
As $m \rightarrow \infty$ from (\ref{c32}), we get
$\lim\limits_{m}\frac{f_{m+1}^{2}}{f_{m+1}f_{m}}.\frac{\lambda_{m}}{\lambda_{m}-\lambda_{m-1}}.a_{nm}y_{m}$ exists $(n \in \mathbb{N})$
and since $y \in c\smallsetminus c_{0},$ we have
$\lim\limits_{m}\frac{f_{m+1}^{2}}{f_{m+1}f_{m}}.\frac{\lambda_{m}}{\lambda_{m}-\lambda_{m-1}}.a_{nm}$ exists; $(n \in \mathbb{N})$ which shows the necessity of (\ref{c28}) holds, where $l=\lim\limits_{k}y_{k}.$
\\ 
Again since $Ax \in l_{p}$ and $\bar{A}x\in l_{p}$ implies 
$\left\lbrace a_{n}\right\rbrace \in l_{p} $ by (\ref{c34}).
\par This completes the proof of part $(i).$
\par Part $(ii)$ can be proved in the similar way of that used in the proof of part $(i)$ by using Lemma \ref{l4}. 
\end{proof}
 
\begin{thm}\label{T6}
(i) Let $1 \leq p < \infty .$ Then $A \in \left( c^{\lambda}_{0}\left( \hat{F}\right): l_{p} \right) $ if and only if (\ref{c25}) and (\ref{c26}) hold, and 
\begin{equation}\label{c37}
\sum\limits_{j=k}^{\infty} a_{nj}  \mbox{~exists~} (n,k \in \mathbb{N}),
\end{equation}
\begin{equation}\label{c38}
\left\lbrace \frac{f_{k+1}^{2}}{f_{k}f_{k+1}}. \frac{\lambda_{k}}{\lambda_{k}-\lambda_{k-1}}.a_{nk}\right\rbrace _{k=0}^{\infty} \in l_{\infty} \ (n \in \mathbb{N}).
\end{equation}
(ii) $A \in \left( c^{\lambda}_{0}\left( \hat{F}\right): l_{\infty} \right)$ if and only if (\ref{c30}),(\ref{c37}) and (\ref{c38}) hold
\end{thm}
\begin{proof}
The proof is  similar to the proof of Theorem \ref{T2}.
\end{proof}
\begin{thm}\label{T5}
$A \in \left( c^{\lambda}(\hat{F}): c\right) $ if and only if 
(\ref{c27}),(\ref{c28}),(\ref{c30}) hold and
\begin{equation}\label{c49}
\lim\limits_{n}a_{n}=a,
\end{equation}
\begin{equation}\label{c50}
\lim\limits_{n}\bar{a}_{nk}=\alpha_{k} \:(k \in \mathbb{N})
\end{equation}
and
\begin{equation}\label{c51}
\lim\limits_{n}\sum\limits_{k}\bar{a}_{nk}=\alpha 
\end{equation}
\end{thm} 
 
\begin{proof}
Suppose that $A$ satisfies the conditions (\ref{c27}),(\ref{c28}),(\ref{c30}),(\ref{c49}),(\ref{c50}), (\ref{c51}) and  take any $x \in c^{\lambda}(\hat{F}).$ Since (\ref{c30}) implies (\ref{c26}) we  have by Theorem \ref{T1} that
 $\left\lbrace a_{nk}\right\rbrace_{k \in \mathbb{N}}  \in  \left\lbrace c^{\lambda}(\hat{F})\right\rbrace^{\beta} $ for all $n \in \mathbb{N}$ and hence $Ax$ exists. Also from (\ref{c30}) and (\ref{c50}) we have 
 $$ \sum\limits_{j=0}^{k}\left| \alpha_{j}\right| \leq 
 \sup _{n}\sum\limits_{j}\left| \bar{a}_{nj}\right| < \infty \mbox{~for~all~} k \in \mathbb{N}. $$  This implies that $(\alpha_{k})\in l_{1}$ and hence the series $\sum\limits_{k}\alpha_{k}(y_{k}-l) $
 converges where $y=(y_{k})\in c$ is the sequence connected with $x=(x_{k})$ by $y=\bar{F}(x)$ such that $y_{k}\rightarrow l$ as
 $k \rightarrow \infty.$ By combining Lemma \ref{l3} with the conditions (\ref{c30}),(\ref{c50}) and (\ref{c51}) that the matrix $A=(\bar{a}_{nk})$ is in the class $(c :c).$
 Now, applying the same method used in Theorem \ref{T2} we obtain that relation (\ref{c33}) holds and from which we obtain
 \begin{equation}\label{c52}
 \sum\limits_{k} a_{nk}x_{k}
 =\sum\limits_{k} \bar{a}_{nk}\left( y_{k}-l\right) 
 +l \sum\limits_{k} \bar{a}_{nk}+la_{n} \: (n \in \mathbb{N}).
 \end{equation}
 As $n \rightarrow \infty$ from (\ref{c52}) we have,
 $$ A_{n}(x) \rightarrow \sum\limits_{k} \bar{a}_{nk}\left( y_{k}-l\right)+l\alpha+la , $$
 which shows $A(x) \in c$ and hence $A \in \left(c^{\lambda}(\hat{F}) : c \right) .$
 \par Conversely suppose $A \in \left(c^{\lambda}(\hat{F}) : c \right).$ Since $c \subset l_{\infty},$ we have 
 $A \in \left(c^{\lambda}(\hat{F}) : l_{\infty} \right).$ This leads us with Theorem \ref{T2} to the necessity conditions (\ref{c27}),(\ref{c28}) and (\ref{c30}). Consider the sequence 
 $b^{(k)}=\left\lbrace b^{(k)}_{n}\right\rbrace_{n \in \mathbb{N}}\in c^{\lambda}(\hat{F}) $ defined by (\ref{c5}) for all fixed $k \in \mathbb{N}.$ Then we have 
 $Ab^{(k)}=\left\lbrace \bar{a}_{nk}\right\rbrace_{n \in \mathbb{N}} $ and hence $\left\lbrace \bar{a}_{nk}\right\rbrace_{n \in \mathbb{N}}\in c$ for all $k \in \mathbb{N}$ which gives condition (\ref{c50}). Let $z=\sum\limits_{k} b^{(k)}.$ Then the linear transformation $T :c^{\lambda}(\hat{F}) \rightarrow c $ where $T \equiv\bar{F}$ is continuous and we obtain that $\bar{F}_{n}(z)=\sum\limits_{k} \bar{F}_{n}\left( b^{(k)}\right)=1,\: (n \in \mathbb{N}) $ which gives $\bar{F}(z)=e \in c$ and hence $z\in c^{\lambda}(\hat{F}).$
 \par Since $c^{\lambda}(\hat{F})$ and $c$ are $BK-$spaces, therefore Lemma \ref{l5} implies the continuity of the matrix mapping $A:c^{\lambda}(\hat{F}) \rightarrow c .$ Thus we have for every $n \in \mathbb{N},$ that 
 $$ A_{n}(z)=\sum\limits_{k} A_{n}\left( b^{(k)}\right) =\sum\limits_{k} \bar{a}_{nk}.$$
 This shows the necessity of (\ref{c51}).
 \par Now, it follows that (\ref{c30}), (\ref{c50}) and (\ref{c51})
 with Lemma \ref{l3} that $\bar{A}=(\bar{a}_{nk})\in (c :c).$ This leads us with (\ref{c27}) and (\ref{c28}) to the consequence that the relation (\ref{c34}) holds for all $x \in c^{\lambda}(\hat{F})$ and $y \in c$ which are connected by
 $y=\bar{F}(x)$ such that $y_{k} \rightarrow l$ as $k \rightarrow \infty.$
 \par Since  $Ax \in c$ and $\bar{A}y \in c, $  the necessity of (\ref{c49}) is obtained by (\ref{c34}).
\end{proof}
 
 \begin{thm}
 $A \in \left( c^{\lambda}(\hat{F}): c_{0}\right) $ if and only if 
 (\ref{c27}),(\ref{c28}),(\ref{c30}) hold and
 \begin{equation}\label{c53}
 \lim\limits_{n}a_{n}=0,
 \end{equation}
 \begin{equation}\label{c54}
 \lim\limits_{n}\bar{a}_{nk}=0 \:(k \in \mathbb{N})
 \end{equation}
 and
 \begin{equation}\label{c56}
 \lim\limits_{n}\sum\limits_{k}\bar{a}_{nk}=0 
 \end{equation}
 \end{thm} 
 \begin{proof}
 It  can be proved in the similar way as Theorem \ref{T5} with Lemma \ref{l7}.
 \end{proof}
 
\begin{thm}\label{T7}
$A \in \left( c_{0}^{\lambda}(\hat{F}): c\right) $ if and only if 
(\ref{c30}),(\ref{c37}),(\ref{c38}) and (\ref{c50}) hold.
\end{thm}
\begin{proof}
This result can be proved by using Lemma \ref{l2}, Theorem \ref{T1}
and Theorem \ref{T6}.
\end{proof} 

\begin{thm}
$A \in \left( c_{0}^{\lambda}(\hat{F}): c_{0}\right) $ if and only if 
(\ref{c30}),(\ref{c37}),(\ref{c38}) and (\ref{c54}) hold.
\end{thm}
\begin{proof}
This result can be proved by using Lemma \ref{l8}, Theorem \ref{T1}
and Theorem \ref{T7}.
\end{proof} 
\begin{lem} \cite{Altay,16} \label{l10}
Let $X$ and $Y$ be any two sequence spaces, $A$ is an infinite matrix and $B$ be a triangle. Then $A \in \left( X : Y_{B}\right) $
if and only if $BA \in \left( X :Y\right). $
\begin{cor}
Let $A=(a_{nk})$ be an infinite matrix and define the matrix $C=(c_{nk})$ by
$$ c_{nk}=\frac{1}{\lambda_{n}}\sum\limits_{i=0}^{n}\left(\lambda_{i}-\lambda_{i-1} \right)\left( \frac{f_{i}}{f_{i+1}}a_{ik}-\frac{f_{i+1}}{f_{i}}a_{i-1,k}\right); \; (n,k \in \mathbb{N}).  $$
By applying Lemma \ref{l10} we get,
 $A$ belongs to any one of the classes
 $\left( c_{0}:c_{0}^{\lambda}(\hat{F}) \right),$
 $\left( c:c_{0}^{\lambda}(\hat{F}) \right),$
 $\left( l_{p}:c_{0}^{\lambda}(\hat{F}) \right),$
 $\left( c_{0}:c^{\lambda}(\hat{F}) \right),$
 $\left( c:c^{\lambda}(\hat{F}) \right),$
 $\left( l_{p}:c^{\lambda}(\hat{F}) \right)$
 if and only if the matrix $C$ belongs respectively to the classes
 $\left( c_{0}:c_{0} \right),$
  $\left( c:c_{0} \right),$
  $\left( l_{p}:c_{0} \right),$
  $\left( c_{0}:c \right),$
  $\left( c: c\right),$
  $\left( l_{p}:c \right),$ where $1\leq p \leq \infty.$
 
\end{cor}
\end{lem} 
\section{ Compact operators on the spaces $c_{0}^{\lambda}(\hat{F})$ and $c^{\lambda}(\hat{F})$} 
  
\par In this section, we establish some estimates for the operator norms and the Hausdorff measures of noncompactness of certain matrix operators on the spaces 
$c_{0}^{\lambda}(\hat{F})$ and $c^{\lambda}(\hat{F}).$ Further, by using the Hausdorff measure of noncompactness, we characterized some classes of compact operators on these spaces.\\

 For our investigations we need the following results.
\begin{thm} \cite{cp3, cp6} \label{o1}
Let $X$ and $Y$ be FK spaces. Then $(X,Y) \subset B(X,Y),$
that is, every $A \in (X,Y)$ defines a linear operator $L_{A} \in B(X,Y)$ where $L_{A}(x)=A(x)$ and $x \in X.$
\end{thm}
 
\begin{thm} \cite{cp4}\label{o2}
Let $X \supset \phi$ and $Y$ be BK spaces. Then  $A \in \left( X, l_{\infty}\right) $if and only if 
 $\parallel  A \parallel_{X}^{*} 
=\sup_{n} \parallel  A_{n} \parallel_{X}^{*}< \infty.$\\
Furthermore, if  $A \in \left( X, l_{\infty}\right) $ then it follows that 
$\parallel L_{A}\parallel=\parallel A\parallel_{X}^{*}.$
\end{thm} 

\begin{thm} \cite{cp2} \label{o5}
Let $X$ be a BK space. Then $A \in \left( X,l_{1}\right) $ if and only if \\
$\parallel A \parallel^{*}_{X,1}
=\sup_{N \subset \mathbb{N}} \parallel\left( \sum\limits_{n \in N} a_{nk}\right) _{k=0}^{\infty} \parallel_{X}^{*} < \infty,$ where $N$ is finite.\\
Moreover, if $A \in \left( X, l_{1}\right) $ then
$\parallel A \parallel^{*}_{X,1}\leq \parallel L_{A}\parallel\leq 4.\parallel A \parallel^{*}_{X,1}.$
\end{thm}
\par Throughout, let $T=(t_{nk})_{n,k=0}^{\infty}$ be a triangle, that is $t_{nk}=0$ for $k>n$ and $t_{nn}\neq 0 , (n=0,1,...)$ and $S$ its inverse. The following results are known.
\begin{thm} \cite{cp1,cp6} \label{o3}
Let $\left( X, \parallel . \parallel \right) $ be a BK space. Then $X_{T}$ is a BK space with 
$\parallel . \parallel_{T}= \parallel T(.)\parallel.$
\end{thm}
\begin{rem} \cite{cp1}
The matrix domain $X_{T}$ of a normed sequence space $X$ has basis if and only if $X$ has a basis.
\end{rem}

\begin{thm} \cite{cp1}\label{T3}
Let $X$ be a BK space with AK and $R=S^{t},$ the transpose of $S.$
If $a \in \left( X_{T}\right)^{\beta} $ then 
$\sum\limits_{k=0}^{\infty}a_{k}x_{k}=\sum\limits_{k=0}^{\infty}R_{k}(a)T_{k}(x)$ for all $x \in X_{T}.$
\end{thm}

\begin{rem} \cite{cp1}
The conclusion of Theorem \ref{T3} holds for $X=c$ and $X=l_{\infty}.$
\end{rem}

\begin{thm} \cite{cp3}
Let $X$ and $Y$ be Banach spaces, $S_{X}=\left\lbrace x \in X : \parallel x \parallel=1 \right\rbrace ,$
$K_{X}=\left\lbrace x \in X : \parallel x \parallel \leq 1 \right\rbrace$ and $A\in B(X,Y).$ Then the Hausdorff measure of noncompactness of a compact operator $A,$ denoted by $\parallel A \parallel_{\chi},$ is given by 
$\parallel A \parallel_{\chi}=\chi\left( AK \right)=\chi\left( AS \right). $
\end{thm}
\par Furthermore, $A$ is compact if and only if $\parallel A \parallel_{\chi}=0$ (see \cite{cp3}). The Hausdorff measure of noncompactness satisfies the inequality 
 $\parallel A \parallel_{\chi} \leq  \parallel A \parallel$ 
 (see \cite{cp3}).

\begin{thm} \cite{cp3} \label{o7}  
Let $X$ be a Banach space with Schauder  basis $\left\lbrace
e_{1}, e_{2}, ... \right\rbrace, $ $Q$ be a bounded subset of $X,$
and $P_{n}: X\rightarrow X$ be the projector onto the linear span of $\left\lbrace
e_{1}, e_{2}, ...,e_{n} \right\rbrace.$ Then 
$$\frac{1}{a} \limsup\limits_{n \rightarrow \infty} \left( \sup_{x \in Q}\parallel (I-P_{n})x\parallel\right)  \leq \chi(Q)
\leq \limsup\limits_{n \rightarrow \infty} \left( \sup_{x \in Q}\parallel (I-P_{n})x\parallel\right),$$
where $$ a=\limsup\limits_{n \rightarrow \infty} \parallel I-P_{n}\parallel.$$
\end{thm} 

\begin{thm} \cite{cp5} \label{o6}
Let $Q$ be a bounded subset of a normed space $X,$ where $X$ is $l_{p}$ for $1 \leq p < \infty$ or $c_{0}.$ If  $P_{n}: X\rightarrow X$ is an operator defined by $P_{n}(x)=\left( x_{0},x_{1},...,x_{n}, 0, 0, ... \right), $ then 
$$\chi(Q)= \lim \limits_{n \rightarrow \infty} \left( \sup_{x \in Q}\parallel (I-P_{n})x\parallel\right).$$
\end{thm}

\begin{thm} \cite{cp7}
Let $X$ be a normed sequence space and $\chi_{T}$ and $\chi$ denote the Hausdorff measures of noncompactness on $M_{X_{T}}$ and $M_{X},$ the collection of all bounded sets in $X_{T}$ and $X,$ respectively. Then $\chi_{T}(Q)=\chi(T(Q))$ for all $Q \in M_{X_{T}}.$
\end{thm}

\begin{lem} \cite{cp3} \label{l9}
 Let $X$ denote any of the spaces $c_{0},c$ or $l_{\infty}.$ Then $x^{\beta}=l_{1}$ and $\parallel a \parallel_{X}^{*}
=\parallel a\parallel_{l_{1}}$ for all  $a \in l_{1}.$
\end{lem}

\begin{thm}
(a) Let $A \in \left( c^{\lambda}(\hat{F}), l_{\infty}\right), $ the matrix $\bar{A}=\left( \bar{a}_{nk}\right)_{n,k=0}^{\infty} $
be defined by
 \[\bar{a}_{nk}
   =\lambda_{k}\left[ \frac{a_{nk}}{\lambda_{k}-\lambda_{k-1}}. \frac{f_{k+1}^{2}}{f_{k}f_{k+1}}+
   \left( \frac{1}{\lambda_{k}-\lambda_{k-1}}. \frac{1}{f_{k}f_{k+1}}-
   \frac{1}{\lambda_{k+1}-\lambda_{k}}. \frac{1}{f_{k+1}f_{k+2}}
   \right)\sum\limits_{j=k+1}^{\infty}f_{j+1}^{2}a_{nj}
   \right]\]
   for all $n,k=0,1...$ and 
   $\parallel A\parallel_{\left( c^{\lambda}(\hat{F}), l_{\infty}\right) }=\sup_{n}\left( \sum\limits_{k=0}^{\infty}\left| \bar{a}_{nk}\right| \right) .$\\ Then
   \begin{equation}
   \parallel L_{A}\parallel=\parallel A\parallel_{\left( c^{\lambda}(\hat{F}), l_{\infty}\right) }.
   \end{equation}
  (b) Let  $A \in \left( c^{\lambda}(\hat{F}), l_{1}\right)  $ and 
   $\parallel A\parallel_{\left( c^{\lambda}(\hat{F}), l_{1}\right) }=\sup_{N \subset \mathbb{N}}\left( \sum\limits_{k}\left| \sum \limits_{n\in N}\bar{a}_{nk}\right| \right),$ where $N$ is finite.\\
   Then  
    \begin{equation}
    \parallel A\parallel_{\left( c^{\lambda}(\hat{F}), l_{1}\right) } \leq
      \parallel L_{A}\parallel \leq 4.\parallel A\parallel_{\left( c^{\lambda}(\hat{F}), l_{1}\right) }.
      \end{equation}
\end{thm} 
\begin{proof}
(a) We assume that 
$A \in \left( c^{\lambda}(\hat{F}), l_{\infty}\right) ,$ then we have $A_{n} \in \left( c^{\lambda}(\hat{F}) \right)^{\beta} $ for all $n=0,1,2,...,$ and it follows from the Theorem \ref{T3} that
\begin{equation}\label{c39}
A_{n}(x)=\sum\limits_{k=0}^{\infty}a_{nk}x_{k}
=\sum\limits_{k=0}^{\infty}R_{k}(A_{n})T_{k}(x)
\end{equation}
for all $x \in c^{\lambda}(\hat{F})$ and $n=0,1,2,...,$ where
\begin{equation}\label{c40}
R_{k}(A_{n})=\sum\limits_{j=0}^{\infty}r_{kj}a_{nj}
=\sum\limits_{j=0}^{\infty}s_{jk}a_{nj}.
\end{equation}
Here $T=\bar{F}$ and $S=\bar{F}^{-1}.$ Therefore we have

\[R_{k}(A_{n})=
\lambda_{k}\left[ \frac{a_{nk}}{\lambda_{k}-\lambda_{k-1}}. \frac{f_{k+1}^{2}}{f_{k}f_{k+1}}+
   \left( \frac{1}{\lambda_{k}-\lambda_{k-1}}. \frac{1}{f_{k}f_{k+1}}-
   \frac{1}{\lambda_{k+1}-\lambda_{k}}. \frac{1}{f_{k+1}f_{k+2}}
   \right)\sum\limits_{j=k+1}^{\infty}f_{j+1}^{2}a_{nj}
   \right]\]
   i.e.
\begin{equation}\label{c41}
R_{k}(A_{n})=\bar{a}_{nk}
\end{equation}
for all $n$ and $k.$
Since $c^{\lambda}(\hat{F})$ is a BK space, Theorem \ref{o2} gives
\begin{equation}\label{c42}
\parallel A \parallel^{*}_{c^{\lambda}(\hat{F})}=\sup_{n} \parallel A_{n}\parallel^{*}_{c^{\lambda}(\hat{F})}
=\parallel L_{A}\parallel.
\end{equation}
We also have $x \in S_{c^{\lambda}(\hat{F})}$ if and only if $y=\bar{F}(x)\in S_{c}$ by Theorem \ref{o3} and conclude from (\ref{c39}), (\ref{c41}) and the definition of the norms 
$\parallel .\parallel^{*}_{c^{\lambda}(\hat{F})}$ and 
$\parallel .\parallel^{*}_{c},$
\begin{equation}\label{o4}
\parallel A_{n}\parallel^{*}_{c^{\lambda}(\hat{F})}
=\sup\left\lbrace \left| A_{n}(x)\right|: x\in S_{c^{\lambda}(\hat{F})} \right\rbrace 
=\sup\left\lbrace \left| \bar{A}_{n}(y)\right|: y\in S_{c} \right\rbrace 
=\parallel \bar{A}_{n}\parallel^{*}_{c}
\end{equation}
for all $n=0,1,2,....$
By using the Lemma \ref{l9}, we have
\begin{equation}
\parallel L_{A} \parallel= \parallel A \parallel^{*}_{c^{\lambda}(\hat{F})}
=\sup_{n} \parallel A_{n}\parallel^{*}_{c^{\lambda}(\hat{F})}
=\sup_{n}\parallel \bar{A}_{n}\parallel^{*}_{c}
=\sup_{n}\left( \sum\limits_{k=0}^{\infty} \left|\bar{a}_{nk} \right| \right) 
=\parallel A \parallel_{(c^{\lambda}(\hat{F}),l_{\infty})}. 
\end{equation}
\par Part(b) is proved in exactly the same way as part(a), we apply Theorem \ref{o5} instead of Theorem \ref{o2}.
\end{proof}

\begin{thm}\label{T4}
Let $A$ be an infinite matrix and put
$\parallel A \parallel^{(m)}_{\left(c^{\lambda}(\hat{F}), l_{\infty}\right) }
=\sup_{n>m} \sum\limits_{k=0}^{\infty}\left| \bar{a}_{nk} \right|.$
\begin{enumerate}
\item[(a)] If $A \in \left( c^{\lambda}(\hat{F}), c_{0}\right), $ then
\begin{equation}\label{c43}
\parallel L_{A}\parallel_{\chi}=\lim\limits_{m \rightarrow \infty}
\parallel A \parallel^{(m)}_{\left( c^{\lambda}(\hat{F}), l_{\infty}\right)}.
\end{equation}
\item[(b)] If $A \in \left( c^{\lambda}(\hat{F}), c\right), $ then
 \begin{equation}\label{c44}
\frac{1}{2}.\lim\limits_{m \rightarrow \infty}
\parallel A \parallel^{(m)}_{\left( c^{\lambda}(\hat{F}), l_{\infty}\right)}\leq \parallel L_{A}\parallel_{\chi} \leq \lim\limits_{m \rightarrow \infty}
\parallel A \parallel^{(m)}_{\left( c^{\lambda}(\hat{F}), l_{\infty}\right)}.
\end{equation}
\item[(c)] If $A \in \left( c^{\lambda}(\hat{F}), l_{\infty}\right), $ then
 \begin{equation}\label{c45}
0\leq \parallel L_{A}\parallel_{\chi} \leq \lim\limits_{m \rightarrow \infty}
\parallel A \parallel^{(m)}_{\left( c^{\lambda}(\hat{F}), l_{\infty}\right)}.
\end{equation}
\end{enumerate}

\end{thm}
\begin{proof}
Let us assume that limits in (\ref{c43}) to (\ref{c44}) exist. We write $K=\left\lbrace x \in  c^{\lambda}(\hat{F}): \parallel x \parallel \leq 1\right\rbrace $
($a$) Applying Theorem \ref{o6}, we have
\begin{equation}\label{c46}
\parallel L_{A}\parallel_{\chi}=\chi(AK)
=\lim\limits_{m\rightarrow \infty}\left\lbrace  \sup_{x \in K}\parallel \left( I-P_{m}\right)Ax \parallel \right\rbrace,
\end{equation}
where $P_{m}: c_{0}\rightarrow c_{0} \:\left( m=0,1,2,...\right) $
is the projector such that $P_{m}=\left( x_{0}, x_{1},...,x_{m},0,0,...\right) $ for $x=(x_{k})\in c_{0}.$ It is known that $\parallel I-P_{m}\parallel=1$ for all $m.$ Let $A_{m}=\left( \acute{a}_{nk}\right) $ be the infinite matrix with
 $$\acute{a}_{nk}
 = \left\{
          \begin{array}{ll}
                     0
                    & , 0 \leq n \leq m \\
                    \bar{a}_{nk}
                    & , m<n.
         \end{array}
   \right. $$ 
We have
\begin{equation}\label{c47}
\sup_{x \in K}\parallel \left( I-P_{m}\right)Ax \parallel
= \parallel L_{A_{(m)}}\parallel 
= \parallel A_{(m)}\parallel_{\left(c^{\lambda}(\hat{F}), l_{\infty} \right) } 
= \parallel A\parallel_{\left(c^{\lambda}(\hat{F}), l_{\infty} \right) }^{(m)},
\end{equation}   
and we obtain (\ref{c43}) from (\ref{c46}) and (\ref{c47}). \\
(b) This proof can be done exactly in the same way as in (a).
If $P_{m}: c \rightarrow c \: (m=0,1,2,...)$ is the projector such that $P_{m}(x)=le+\sum\limits_{k=0}^{m}(x_{k}-l)e^{(k)}$ then $\parallel I-P_{m} \parallel =2$ for all $m.$
By applying Theorem \ref{o7} we can prove part (b).\\
($c$) Let $P_{m}: l_{\infty}\rightarrow l_{\infty} \: (m=0,1,2,...)$ by $P_{m}(x)=\left( x_{0},x_{1}, ..., x_{m},... \right) $ for $x=(x_{k})\in l_{\infty}.$ Since 
$AK \subset P_{m}(AK)+ (I-P_{m})(AK),$ by applying properties of $\chi,$ we obtain\\
$\chi(AK) \leq \chi\left( P_{m}(AK)\right) +\chi\left( (I-P_{m})(AK)\right)=\chi\left( (I-P_{m})(AK)\right)
\leq \sup_{x \in K}\parallel \left( I-P_{m}\right)Ax \parallel$\\
i.e. $
0 \leq\chi(AK)\leq \parallel L_{A_{(m)}}\parallel
 = \lim\limits_{m \rightarrow \infty}
\parallel A \parallel^{(m)}_{\left( c^{\lambda}(\hat{F}), l_{\infty}\right)}. $
\end{proof}   
\begin{cor}
If either $A \in \left( c^{\lambda}(\hat{F}),c\right) $ or 
$A \in \left( c^{\lambda}(\hat{F}),c_{0}\right),$ then $L_{A}$ is compact if and only if 
\begin{equation}\label{c48}
\lim\limits_{m \rightarrow \infty}
\parallel A \parallel^{(m)}_{\left( c^{\lambda}(\hat{F}), l_{\infty}\right)}=0.
\end{equation}
If $A \in \left( c^{\lambda}(\hat{F}),l_{\infty}\right),$ then $L_{A}$ is compact if the condition (\ref{c48}) holds.
\end{cor}
\begin{thm} 
Let $A$ be an infinite matrix and put\\
$$\parallel A \parallel^{(m)}_{\left( c^{\lambda}(\hat{F}), l_{1}\right)}= \sup_{N \subseteq \mathbb{N}\backslash \left\lbrace 0,1,...,m \right\rbrace}\left( \sum\limits_{k} \left| \sum\limits_{n \in N} \bar{a}_{nk}\right| \right) 
,$$ $N$ is finite. Then
$$\parallel L_{A} \parallel_{\chi}=\lim\limits_{m \rightarrow \infty}
\parallel A \parallel^{(m)}_{\left( c^{\lambda}(\hat{F}), l_{1}\right)}.$$
\end{thm}
\begin{proof}
The proof is similar to that of Theorem \ref{T4}.
\end{proof}
\begin{cor}
If $A \in \left( c^{\lambda}(\hat{F}),l_{1}\right),$ then $L_{A}$ is compact if and only if 
$$\lim\limits_{m \rightarrow \infty}
\parallel A \parallel^{(m)}_{\left( c^{\lambda}(\hat{F}), l_{1}\right)}=0.$$
\end{cor}
\begin{cor}
(a) If $A \in \left( c^{\lambda}(\hat{F}), bv \right) $ then $L_{A}$ is compact if and only if
$$\lim\limits_{m \rightarrow \infty} \sup_{N \subseteq \mathbb{N}\backslash \left\lbrace 0,1,...,m \right\rbrace}\left( \sum\limits_{k} \left| \sum\limits_{n \in N}\left(  \bar{a}_{nk}-  \bar{a}_{n-1,k}\right) \right| \right)=0. $$
\end{cor}
\begin{rem}
The previous results would be true for $c_{0}^{\lambda}(\hat{F})$ instead of $c^{\lambda}(\hat{F}).$
\end{rem}

\end{document}